\documentclass[showpacs, aps]{revtex4}

\usepackage{amsmath}
\usepackage{amssymb}
\usepackage{graphicx}
\usepackage{amssymb}
\newcommand{\R}{\mathbb{R}}

\begin{document}

\title{Arrest of Langmuir wave collapse by quantum effects}
\author{
G. Simpson$^1$, C. Sulem$^1$ and P.L. Sulem$^2$}
\affiliation{
$^1$ Department of Mathematics, University of Toronto, 
40 St. George St., Toronto M5S 2E4, Canada\\
$^2$ Universit\'e de Nice-Sophia Antipolis,
CNRS, Observatoire de la C\^ote d'Azur, B.P. 4229, 
06304 Nice Cedex 4, France } 

\begin{abstract}
The arrest of Langmuir-wave collapse by quantum effects, first addressed by Haas and Shukla [Phys. Rev. E  79, 066402
(2009)] using a Rayleigh-Ritz trial-function method is revisited, using rigorous estimates and systematic 
asymptotic expansions. The absence of blow up for the so-called quantum Zakharov equations
is proved in two and three dimensions,  whatever the strength of the quantum effects. 
The time-periodic behavior of the solution for initial
conditions slightly in excess of the singularity threshold for the classical problem is established for various settings
in two space dimensions. The difficulty of developing a consistent perturbative approach in three dimensions is also discussed,
and a semi-phenomenological model is suggested for this case.
\end{abstract}


\maketitle

\section{Introduction}

Special interest was recently devoted to quantum corrections to the Zakharov equations for Langmuir waves 
in a plasma \cite{Zakh72}. First
considered in one space dimension \cite{G05}, the model was then 
extended to  two and three dimensions\cite{HS09}, 
in a formulation retaining magnetic field fluctuations \cite{Kuz74}. 
In a non-dimensional form, the equations that govern the amplitude ${\mathbf E}$  of the
electric field oscillations and the number density  $n$ read
\begin{eqnarray}
&& i\partial_t {\mathbf E} - \alpha\mbox{\boldmath $\nabla$} \times 
( \mbox{\boldmath $\nabla$}\times {\mathbf E})  +  
\mbox{\boldmath $\nabla$}(\mbox{\boldmath $\nabla$} \cdot {\mathbf E})  \nonumber\\
&&\qquad \quad  = n {\mathbf E} + 
\Gamma  \, \mbox{\boldmath $\nabla$}\Delta (\mbox{\boldmath $\nabla$} \cdot {\mathbf E})  \label{ZV1}\\
&&\partial_{tt} n - \Delta  n = \Delta |{\mathbf E}|^2 -\Gamma \Delta^2 n. \label{ZV2}
\end{eqnarray}
In the above equations, the parameter $\alpha$ defined as the square ratio of the light speed and 
the electron Fermi velocity is usually large.
The corresponding term is nevertheless  moderate, because, 
magnetic effects are relatively weak, making ${\mathbf E}$ close to a gradient field. 
In contrast, the coefficient  $\Gamma$ that measures the influence of quantum effects 
is usually very small. We refer to  \cite{HS09} for a discussion of the physical regimes described
by the present model and an estimate of the plasma parameters.
Typically, for a hydrogen plasma, one has $\alpha \approx 8.  10^{23} \,n_0^{-2/3}$ and 
$\Gamma \approx 6. 10^6 \,n_0^{-1/3}$, leading to   $\alpha \approx 8. 10^2$ and $\Gamma \approx 10^{-4}$
in the case of  the rather high equilibrium densities  $n_0= 10^{32} \, m^{-3}$.

Equation (\ref{ZV2}) originates from  the hydrodynamic system
\begin{eqnarray}
&&\partial_t n +  \mbox{\boldmath $\nabla$}\cdot {\mathbf v}=0\\
&&\partial_t  {\mathbf v} = - \mbox{\boldmath $\nabla$} (n + |{\mathbf E}|^2) + 
\Gamma  \mbox{\boldmath $\nabla$} \Delta n,
\end{eqnarray}
governing the ion sound waves.

For  $\Gamma = 0$, it was  shown {\cite  {SuSu79,AA84,  GM94b, GTV97}} that for  $\alpha\geqslant 1$
and ``small enough'' initial conditions, the solution remains smooth for all time.
In two dimensions, the smallness
condition reads $|{\mathbf E}_0|_{L2}2 \le |R|^2_{L2} \approx 1.86$
(where $R$ is the ground state defined in (\ref{NLSgrst})) and is
optimal. In three dimensions, it requires that 
the plasmon number ${\cal N}$ and the  Hamiltonian ${\cal H}$ (defined in (\ref{numb}),(\ref{Ham}))
satisfy ${\cal N}|{\cal H}| <2.6 \ 10^{-4}$ 
together with $|\mbox{\boldmath $\nabla$} {\mathbf E}_0|_L^2 < |{\cal H}|$, conditions that are probably
much too strict.
% A more refined analysis requiring minimal assumptions 
%on the initial conditions is given in \cite{BC96}.

Although not rigorously proved,
the phenomenon of wave collapse is expected for $\Gamma =0$, when the initial conditions are large enough. 
The existence of a finite-time blow up is indeed suspected when the Hamiltonian is negative,
on the basis of numerical simulations and heuristic  arguments (see \cite{SuSu99} for review).

The question then arises of the possible arrest of collapse by 
quantum corrections. This issue was addressed in \cite{HS09} by implementing an
approach based on the Rayleigh-Ritz trial function method, in regimes where the quantum
Zakharov equations can  be reduced to a vector nonlinear Schr\"odinger equation. The latter
results from the  assumption that the density  is slaved to the electric field  oscillations
(adiabatic approximation)
\begin{equation}
- \Delta  n = \Delta |{\mathbf E}|^2 -\Gamma \Delta^2 n
\end{equation}
and, because of the smallness of $\Gamma$, can  be expressed  to leading order as
$\displaystyle{n= - |{\mathbf E}|^2 - \Gamma \Delta |{\mathbf E}|^2}$.
Although this approach led to interesting conclusions such as the arrest of collapse
by quantum effects and its replacement by a time-periodic solution,
it nevertheless involves possibly questionable assumptions.
For small enough $\Gamma$, quantum effects  become relevant 
only very close to the singularity when the adiabatic approximation, even if 
valid at early times, hardly holds.
The rates of blow-up of the solutions of the Zakharov equations (with $\Gamma =0$) are indeed such that
all the terms in eq. (\ref{ZV2}) have the same magnitude  in two dimensions, while 
$\partial_{tt} n \gg \Delta n$ in three dimensions (supersonic collapse).
Furthermore, the  
Rayleigh-Ritz trial function method used to  reduce the problem to the evolution of 
a few scaling coefficients, is based on an arbitrary choice the functional form
of the solution.  The aim 
of the present approach is to revisit the issue of collapse arrest by 
quantum corrections, mainly in regimes amenable to a systematic analysis. 
The paper is organized as follows.  Section II provides a rigorous proof 
of the arrest of Langmuir
collapse by arbitrarily small quantum effects,  in the general framework of the 
Zakharov equations (\ref{ZV1},\ref{ZV2}), both in two and three space dimensions. 
Section III reviews the electrostatic approximation that is valid when the plasma
is not too hot, as well as  the so called scalar model \cite{DKO84} that in the case of
rotational symmetry does not prescribe a zero electric field at the symmetry center 
and a ring (2D) or shell (3D) profile for the electric field intensity.
Section IV provides an asymptotic analysis of the 
dynamics in the presence of weak quantum  effects for the scalar model in two space dimensions. 
This analysis is extended to the two-dimensional electrostatic equations with rotational
symmetry in Section V. The difficulties of the 
three-dimensional problem are discussed in Section VI, where a phenomenological model
based on a heuristic extension of perturbative calculations is presented.
Section VII briefly summarizes our conclusions.

\section{Arrest of collapse by quantum effects}

In this section, we present a rigorous proof of the absence of wave collapse for the general
Zakharov equations with quantum effects (\ref{ZV1}, \ref{ZV2}) in space  dimension $d=2$
and $3$. For this purpose, we first define the  $L^p(\R^d)$  and  Sobolev spaces   
 $H^s(\R^d)$ as the spaces of
scalar,  vector or tensor functions  respectively equipped with the norms \cite{Adams78} 
\begin{eqnarray}
&&|f|_{L^p}= \Big(\int |f|^p d{\mathbf x}\Big )^{1/p}\\
&&|f|_{H^s}= \Big(\int (1+ k^{2})^s|{\widehat f(k)}|^2 d{\mathbf x}\Big )^{1/2},
\end{eqnarray} 
where $\displaystyle{\widehat f(k)}$ denotes the spatial Fourier transform of the function $f$. 
We will also  use of a special case of the  Gagliardo-Nirenberg inequality  (see e.g.
\cite{C03} for the general formula), in the form:
for $\displaystyle{f \in \{L^q(\R^d)},   \mbox{\boldmath $\nabla$}f \in L^2(\R^d)\}$, one has
\begin{equation}
|f|_{L^p} \leqslant K\, |\mbox{\boldmath $\nabla$}f|_{L^2}^\theta\,|f|_{L^q}^{1-\theta},
\end{equation}
where $K$ is a positive constant whose optimal value is given in \cite{W83} in two space dimensions
and in a more general setting in \cite{A08}. 
In dimensions $d=1$ and $2$, $1 < q < p$, 
while for $d > 2$, $\displaystyle{1 < q < p < \frac{2d}{d-2}}$. In all the cases, 
$\theta = \displaystyle {\frac{2 d (p-q)}{p[2d-q(d-2)]}}$.

Among the conserved quantities of the  quantum Zakharov equations (\ref{ZV1},\ref{ZV2}), 
the number of plasmons and the Hamiltonian play an important role 
in the regularity properties of the solution. They read \cite{HS09}
\begin{eqnarray}
&&{\cal N} =  \int |{\mathbf E}|^2 d{\mathbf x} \label{numb}\\
&&{\cal H} = \int \Big \{  |\mbox{\boldmath $\nabla$}\cdot {\mathbf E}|^2 +
\alpha  |\mbox{\boldmath $\nabla$}\times {\mathbf E}|^2 
+ \frac{1}{2} n^2 + \frac{1}{2}|{\mathbf v}|^2\nonumber \\
&&\qquad + n |{\mathbf E}|^2 +\Gamma |\mbox{\boldmath $\nabla$}
(\mbox{\boldmath $\nabla$}\cdot {\mathbf E})|^2  + \frac{\Gamma}{2}
|\mbox{\boldmath $\nabla$}n|^2 \Big \} d{\mathbf x} \label{Ham}. 
\end{eqnarray}

Among all the terms arising in the Hamiltonian, only $ \int n|{\mathbf E}|^2
d{\mathbf x}$ could be non-positive and thus needs to be estimated.
One has (denoting by $C$ different constants)
\begin{eqnarray}
&&\int n\,|{\mathbf E}|^2 d{\mathbf x} \leqslant |n|_{L^4} \,\big | |{\mathbf E}|^2  \big |_{L^{4/3}}
=  |n|_{L^4} \,|{\mathbf E}|^2_{L^{8/3}} \nonumber \\
&&\qquad  \leqslant C |n|_{H^1} \,|{\mathbf E}|^2_{L^{8/3}}\leqslant
\frac{\Gamma}{4} |n|^2_{H^1} + \frac{C}{\Gamma}|{\mathbf E}|^4_{L^{8/3}}, 
\end{eqnarray}
where the successive inequalities result from the H\"older,  Gagliardo-Nirenberg
and  Young inequalities. Using again the  Gagliardo-Nirenberg inequality, we write
\begin{equation}
|{\mathbf E}|_{L^{8/3}}  \leqslant C |\mbox{\boldmath $\nabla$} {\mathbf E}|^{d/8}_{L^2}
\,|{\mathbf E}|^{1-d/8}_{L^2}. 
\end{equation}
It follows that 
\begin{equation}
\int n\,|{\mathbf E}|^2 d{\mathbf x} \leqslant \frac{\Gamma}{4} |n|_{H^1}^2 + \frac{C}{\Gamma}
{\cal N}^{2-d/4}|\mbox{\boldmath $\nabla$} {\mathbf E}|^{d/2}_{L^2}. \label{integr}
\end{equation}

It is then convenient to rewrite the Hamiltonian in the form 
\begin{eqnarray}
&&{\cal H} =    |\mbox{\boldmath $\nabla$} {\mathbf E}|_{L^2}^2
+(\alpha-1)  |\mbox{\boldmath $\nabla$}\times {\mathbf E}|_{L^2}^2 
+\Gamma |\mbox{\boldmath $\nabla$}
(\mbox{\boldmath $\nabla$}\cdot {\mathbf E})|_{L^2}^2  + \nonumber \\
&&\qquad\frac{\Gamma}{2} |n|_{H^1}^2 + 
\Big (\frac{1}{2}-\frac{\Gamma}{2}\Big)|n|_{L^2}^2
+ \frac{1}{2}|{\mathbf v}|_{L^2}^2  +\int  n |{\mathbf E}|^2
d{\mathbf x}. \nonumber \\
&&
\end{eqnarray}
Using eq. (\ref{integr}), one gets the upper bound 
\begin{eqnarray}
&& |\mbox{\boldmath $\nabla$} {\mathbf E}|_{L^2}^2 + 
(\alpha-1)  |\mbox{\boldmath $\nabla$}\times {\mathbf E}|_{L^2}^2 
+\Gamma |\mbox{\boldmath $\nabla$}
(\mbox{\boldmath $\nabla$}\cdot {\mathbf E})|_{L^2}^2  \nonumber \\
&&\qquad +\frac{\Gamma}{2} |n|_{H^1}^2 +
\Big (\frac{1}{2}-\frac{\Gamma}{2}\Big)|n|_{L^2}^2
+ \frac{1}{2}|{\mathbf v}|_{L^2}^2 \nonumber \\
&&\qquad \leqslant |{\cal H}|
+\frac{\Gamma}{4} |n|_{H^1}^2 +  \frac{C}{\Gamma}
{\cal N}^{2-d/4}|\mbox{\boldmath $\nabla$} {\mathbf E}|^{d/2}_{L^2}.\label{estim}
\end{eqnarray}
Since $\alpha$ and  $\Gamma$ are respectively larger and smaller 
than 1, one in particular has
\begin{eqnarray}
&& |\mbox{\boldmath $\nabla$} {\mathbf E}|_{L^2}^2
\leqslant |{\cal H}| +   \frac{C}{\Gamma}
{\cal N}^{2-d/4}(|\mbox{\boldmath $\nabla$} {\mathbf E}|_{L^2}^2)^{d/4}.
\end{eqnarray}
which implies that, for any $\Gamma >0$, 
$|\mbox{\boldmath $\nabla$} {\mathbf E}|_{L^2}$
remains uniformly bounded in time.
Equation (\ref{estim}) also provides  an uniform bound for
$|n|_{H^1}$,  $|{\bf v}|_{L^2}$ and 
$|\mbox{\boldmath $\nabla$}(\mbox{\boldmath $\nabla$}\cdot  {\mathbf E})|_{L^2}^2 $.
One then  easily derives, using standard methods,
the existence for all time of a classical solution for the quantum 
Zakharov equations, both in two and three dimensions.

\section{Electrostatic limit and scalar model}

The large value of the parameter $\alpha$ makes the magnetic fluctuations actually subdominant, 
leading to the so called electrostatic approximation \cite{Zakh72}. 
For this purpose, it is convenient to derive from  eq. (\ref{ZV1}) the system
\begin{eqnarray}
&& i\partial_t  (\mbox{\boldmath $\nabla$} \cdot {\mathbf E}) 
+ \Delta (\mbox{\boldmath $\nabla$} \cdot {\mathbf E}) =  n (\mbox{\boldmath $\nabla$} \cdot {\mathbf E})  + 
\mbox{\boldmath $\nabla$} n \cdot {\mathbf E}\nonumber \\
&&\qquad  + \Gamma \Delta^2  (\mbox{\boldmath $\nabla$} \cdot {\mathbf E}) \label{div}\\
&& i \partial_t (\mbox{\boldmath $\nabla$} \times {\mathbf E})   
+ \alpha \Delta (\mbox{\boldmath $\nabla$} \times {\mathbf E}) = 
n (\mbox{\boldmath $\nabla$} \times {\mathbf E}) \nonumber \\
&&\qquad + \mbox{\boldmath $\nabla$} n \times {\mathbf E}. \label{rot}
\end{eqnarray}
Even if the initial electric field is a gradient, $\mbox{\boldmath $\nabla$} \times {\mathbf E}$
is driven by the last term in the right hand side of eq. (\ref{rot}). 
Nevertheless, when $\alpha$ is large,
a stationary-phase argument applied to this equation, easily shows that 
$\mbox{\boldmath $\nabla$} \times {\mathbf E}$ saturates at a level that scales like $1/\alpha$.
Thus,  although small, it contributes to eq. (\ref{ZV1}) but not to the Hamiltonian that,
as $\alpha \to \infty$, has a finite limit obtained by neglecting the term 
involving the coefficient  $\alpha$.

Writing $\displaystyle{{\mathbf E} = 
-\mbox{\boldmath $\nabla$}\psi +  \frac{1}{\alpha}{\mathbf E_1}}$ and substituting in eq. (\ref{div}),
one gets to leading order 
\begin{eqnarray}
&& \Delta (i\partial_t \psi + \Delta \psi) = 
\mbox{\boldmath $\nabla$} \cdot (n \mbox{\boldmath $\nabla$} \psi)
+ \Gamma \Delta^3 \psi, 
\label{zakhpot1}\\
&&\partial_{tt}n - \Delta n = \Delta (|\mbox{\boldmath $\nabla$} \psi|^2) - \Gamma \Delta^2 n.
\label{zakhpot2}
\end{eqnarray}
A rigorous proof of the convergence is given in \cite{Gal00}  in the cases where the solution is globally smooth.
The proximity of a singularity is nevertheless not expected to affect the ordering between the solenoidal and 
gradient components of the electric field. 

The analysis of the solution near collapse is often performed, assuming that the fluctuations involve
rotational symmetry \cite{Zakh72}. In this case, introducing   $E= -\partial\psi /\partial r$, 
eqs. (\ref{zakhpot1})-(\ref{zakhpot2}) reduce to 
%\begin{eqnarray}
%&&i\partial_t E + \partial_r \frac{1}{r^{d-1}}\partial_r r^{d-1}E= nE \noindent \\
%&&\qquad + \Gamma \partial_r\frac{1}{r^{d-1}}\partial_r r^{d-1}\partial_r\frac{1}{r^{d-1}}\partial_r r^{d-1} E \\
%&&\partial_{tt} n - \Delta n=\Delta E - \Gamma \Delta^2 n
%\end{eqnarray}
\begin{eqnarray}
&&i\partial_t E + \Delta_r^{(1)} E= nE + \Gamma  {\Delta_r^{(1)}}^2 E \\
&&\partial_{tt} n - \Delta_r n=\Delta_r |E|^2 - \Gamma \Delta_r^2 n
\end{eqnarray}
where $\Delta_r^{(1)}=   \partial_r r^{-(d-1)}\partial_r r^{d-1}$ 
and  $\Delta_r =r^{-(d-1)}\partial_r r^{d-1}\partial_r$. These equations 
are supplemented with the boundary conditions
$\displaystyle{E(0,t)= \partial_rn(0,t)= E(\infty,t) = n(\infty,t)=0}$.

It was early noticed \cite{DKO84} that the assumption of an electric field vanishing
at the center of symmetry (taken as the origin of coordinates) with a shell profile for 
the intensity profile
is hardly consistent with a realistic model. Relaxing this assumption
of zero electric field at the center
while retaining an isotropic intensity profile of the electric field,
is not possible when the detailed dynamics are retained.
It may thus be suitable \cite{DKO84}
to abandon the vector character of the problem, in order to preserve a 
non-zero electric field at the center of the cavity  as suggested by numerical  simulations of  the vector
Zakharov equation near collapse \cite{PSS91}, while keeping the rotational symmetry 
necessary for implementing an asymptotic analysis in the spirit of \cite{Malkin93,FP99}.
We are thus led to consider the influence of quantum effects in the framework of the 
``scalar model''\cite{DKO84}
\begin{eqnarray}
&& i\partial_t E + \Delta E = n  E + \Gamma \Delta^2 E,  \label{zakhsc1}\\
&&\partial_{tt}n - \Delta n = \Delta  |E|^2 - \Gamma \Delta^2 n. \label{zakhsc2}
\end{eqnarray}
with the condition that $E$ and $n$ vanish at infinity and satisfy 
$\displaystyle{\partial_rE(0,t)= \partial_rn(0,t)=0}$.
Isotropic solutions of these equations are not only  stable
but also attractive near collapse \cite{LPS92}.
It is interesting to notice that direct numerical simulations \cite{PSS91} of the 
collapsing solutions of the vector Zakharov equations
(\ref{ZV1}-\ref{ZV2}) with $\Gamma = 0$ indicate that the anisotropy is 
in general rather  moderate 
and the rates of blow up identical to those  of the scalar model \cite{LPS92}.

For $\Gamma=0$ and in the adiabatic limit where $n= -|E|^2$, 
isotropic solutions of eq. (\ref{zakhsc1}) identifies with
the vortex solutions corresponding to a rotational number $m=1$,
of the (scalar) two-dimensional nonlinear Schr\"odinger equation \cite{FG08}.
Taking the adiabatic limit with $\Gamma \neq 0$,
one gets a non local extension of the cubic nonlinear Schr\"odinger equation with both second and
fourth order dispersions studied in \cite{KS00} and \cite{FIP02}.

\section{Asymptotic behavior of the scalar model in two dimensions} 

In the analysis presented in \cite{HS09},  the space dimension
has no qualitative effect, but this is not the case in the framework of a 
systematic perturbative approach
that, to be fully consistent, requires  a small expansion parameter
usually associated with  closeness to a critical regime. 

\subsection{The classical regime ($\Gamma=0$)}

The two-dimensional regime deserves a special attention because it is amenable to  a detailed
mathematical analysis. It was indeed proved that for two-dimensional smooth 
initial conditions such that the initial density $n_0\in L^2$, $\partial_tn_0\in H^{-1}$ and
the initial electric field $E_0\in  H^1$  obeys  the condition 
$\displaystyle{|E_0|_{L^2} \leqslant |R|_{L^2}}$,
where $R$ is the unique positive solution (ground state) of
\begin{equation}
\Delta R -R + R^3= 0, \label{NLSgrst}
\end{equation}
the solution of the classical scalar model ($\Gamma = 0$) exists for all time in these spaces\cite{GM94b}
and is unique \cite{BC96}. Optimal local existence results in spaces of very weak regularity
appear in \cite{BHHT09}. Further regularity properties of the initial conditions are 
also preserved in time.
The ground state $R$ is radially symmetric \cite{GNN81} and 
obeys the relation
$\int( |\mbox{\boldmath $\nabla$}R|^2 -(1/2) R^4 ) r dr =0$, where the right hand side 
can be viewed as the Hamiltonian of the standing wave solution $e^{it} R(|{\mathbf x}|)$
of the nonlinear Schr\"odinger (NLS) equation $i\psi_t + \Delta \psi + |\psi|^2\psi = 0$.
It is furthermore interesting to notice that, 
unlike the two-dimensional NLS equation, the
Zakharov system does not have blowing up solutions of minimal $L^2$-norm.

Although there is no rigorous proof of existence of a finite-time singularity
for larger initial conditions, 
one has the following result \cite{GM94b} (valid both in two and three dimensions).
Suppose ${\cal H}<0$ and that the solution $(E,n,{\mathbf v})$ is radially symmetric. Then either 
$|E|_{H^1} + |n|_{L^2} + |{\mathbf v}|_{L^2} \to \infty$  as $t \to t_\star$ with $t_\star$ finite, 
or $(E,n, {\mathbf v}  )$ exists for all time and $|E|_{H^1} + |n|_{L^2} + 
|{\mathbf v}|_{L^2} \to \infty$ as $t \to \infty$. 
Numerical simulations clearly indicate the stability
of isotropic solutions and that, near collapse, every solution becomes 
locally isotropic \cite{LPS92}.
Furthermore, in two dimensions, there exist exact self-similar solutions of the classical scalar 
model that blow up in a finite time, of the form
\begin{eqnarray}
&&E(r,t)= \frac{1}{a_0 (t_*-t)} P\Big(\frac{r}{a_0(t_* -t)}
\Big)
e^{i\left (\theta + \frac{1}{a_0^2(t_*-t)}-\frac{r^2}{4(t_*-t)}\right )}\nonumber \\
&& \label{selfsimzak1}\\
&&n(r,t) = \frac{1}{a_0^2(t_*-t)^2} M\Big(\frac{r}{a_0 (t_*-t)}\Big),
\label{selfsimzak2}
\end{eqnarray}
under the condition that $(P,M)$ satisfies the system of ordinary
differential equations in the radial variable
\begin{eqnarray}
&&\Delta P-P-MP =0, \label{Ppeq2} \\
&&a_0^2(\eta^2 M_{\eta \eta}+6\eta M_{\eta}+6M)- \Delta M=\Delta P^2.
\label{Nneq2}
\end{eqnarray}
The parameter $a_0$ entering the self-similar solution is not universal.
One has the following results\cite{GM94a} for existence of solutions to 
(\ref{Ppeq2})--(\ref{Nneq2}).
There exists $a_0^+ >0$ such that  $\forall a_0$ with   $ 0< a_0
< a_0^+$, there is a solution $(P_{a_0}, M_{a_0}) \in
H^1 \times L^2$ of  (\ref{Ppeq2})--(\ref{Nneq2}) with $P_{a_0} >0$. 
Furthermore, 
when $a_0 \to 0$,   $(P_{a_0}, M_{a_0})$ tends to $(R, -R^2)$ in
$H^1 \times L^2$ and for all $c > |R|_{L^2}$, 
there exists  $a_{0c} >0$, such that for all  $a_0$ with  $0< a_0
< a_{0c}$ there is a unique solution  $(P_{a_0}, M_{a_0})$
with  $P_{a_0} >0$ and $|P_{a_0}|_{L^2} < c $.

Numerical simulations of the classical scalar model were performed in \cite{LPS92} 
where it is  shown that for a smooth initial condition
with $|E_0|_{L^2}  >|R|_{L^2}$, the solution
of the initial value problem approaches the self-similar blowing-up solution.  
Furthermore, in a series of simulations with initial conditions such that
$ |E_0|_{L^2}$ approaches
$|R|_{L^2}$ from above, it was observed that the estimated value of the parameter $a_0$ monotonically 
decreases to zero (see table 1 of\cite{LPS92}).
This result is consistent with the existence of a sequence of blowing up solutions 
\begin{eqnarray}
E_{a_0} &=&\frac{1}{1-a_0 t} e^{i\Big (\frac{a_0 |r|^2}
{4(a_0 t -1)} + \frac{t}{1-a_0 t}\Big )} P_{a_0} \Big (
\frac {r}{1 -a_0 t}\Big), \\
n_{a_0} &=&  \frac{1}{(1-a_0 t)^2} M_{a_0} \Big (
\frac {r}{1 -a_0 t}\Big ) 
\end{eqnarray}
obtained by choosing $t_* = -\theta =a_0^{-1}$ in eqs. (\ref{selfsimzak1})-(\ref{selfsimzak2}), with the profiles
obeying  (\ref{Ppeq2})-(\ref{Nneq2}) and thus converging to $(R, -R^2)$ as $a_0\to 0$.
For $a_0=0$, the corresponding solution is smooth, consistent with the regularity of the solutions
such that  $|E|_{L^2}= |R|_{L^2}$ \cite{GM94b}.

The above  observation suggests a pertubative analysis of the influence of
quantum effects for initial conditions such  that
% in the case $\Gamma = 0$, 
$ |E_0|_{L^2}$ is slightly above the threshold for collapse
\footnote{A misprint is to be corrected four lines from
bottom in  the second column of
page 7874 of \cite{LPS92}, which should read $2|R|_{L^2} \approx 2.72450$.}.

\subsection{Influence of quantum effects}

The method for constructing a solution in the presence of weak quantum effects is to assume that these
perturbations induce small corrections of the self-similar profile of the collapsing solution,
but modify the scaling parameter $\lambda$ whose time evolution is prescribed by  the conservation
of the plasmon number and of the Hamiltonian. As noted in \cite{FP99}, an alternative variational approach 
based on the existence of a Lagrangian density is a priori possible, using the actual perturbative 
expansion of the fields. We shall not follow this direction here and concentrate on  direct expansions.

Using a dynamical rescaling transformation, we first define 
$\displaystyle{E(r,t) = \lambda^{-1}U(\xi, \tau)}$,
$\displaystyle{n(r,t) = \lambda^{-2}N(\xi, \tau)}$,
$\displaystyle{n_t(r,t) = \lambda^{-3}W(\xi, \tau)}$ and
$\displaystyle{v(r,t) = \lambda^{-2}V(\xi, \tau)}$,
with $\displaystyle{\xi = r/\lambda}$ and $\tau={\int_0^t \lambda(s)^{-2} ds}$. 
Introducing $a= - \lambda_t= -\lambda^{-2} \lambda_\tau $, the rescaled 
quantities obey 
($\displaystyle{ \Delta = \xi^{-1}\partial_\xi\xi \partial_\xi}\, $)
\begin{eqnarray}
&&i[U_\tau + a\lambda (U+\xi U_\xi)]+\Delta U - NU = \Gamma \lambda^{-2}\Delta^2 U \label{ds1} \\
&&N_\tau + a \lambda(2N+\xi N_\xi) = \lambda W = -\lambda \xi^{-1}\partial_\xi(\xi V)\\
&&V_\tau + a \lambda(2V+\xi V_\xi) +\lambda  N_\xi  \nonumber \\
&&\quad \ = -\lambda \partial_\xi |U|^2
+ \Gamma \lambda^{-1} \partial_\xi \Delta N\\
&&W_\tau + a \lambda(3W+ \xi W_\xi) -\lambda \Delta N  \nonumber \\
&&\quad\  = \lambda \Delta |U|^2 -  
\Gamma \lambda^{-1}  \Delta^2 N.
\end{eqnarray}
It is then convenient to write $\displaystyle{U=e^{i\tau} e^{-ia\lambda\xi^2/4} S}$ and to define
$\displaystyle{b= (a \lambda )_\tau + (a\lambda)^2}$. Equation (\ref{ds1}) is replaced by 
\begin{equation}
iS_\tau -S + \Delta S - NS + b\frac{\xi^2}{4} S =  
\frac{\Gamma}{ \lambda^2} e^{ia\lambda\xi^2/4}   \Delta^2 ( e^{-ia\lambda\xi^2/4}S).
\end{equation}
The assumption is then that we are looking for solutions whose time dependency is only through
the scaling parameter $\lambda$, and thus also through the functions $a$ and $b$. This leads to neglect the 
contributions $S_\tau$, $N_\tau$, $V_\tau$ and $W_\tau$ in the above equations. Combining the resulting 
equations for $N$ and $W$, one gets
\begin{equation}
a^2 {\mathcal L}(N) - \Delta N = \Delta |S|^2 - \frac{\Gamma}{\lambda^2} \Delta^2 N.
\end{equation}
where the operator  ${\mathcal L}$ is defined by 
$\displaystyle{{\mathcal L}= \xi^2\partial_{\xi\xi}+ 6 \xi \partial_\xi + 6}$.

For $\Gamma=0$, the self-similar solution is reached as $t$ 
approaches the singularity time $t_\star$ (i.e. $\tau \to \infty$) and $a$ tends to a constant 
$a_0$, which leads to the conclusion that in this regime
$\lambda(t) = a_0 (t_\star -t)$. Furthermore, as mentioned in Section IV A, $a_0$ tends to zero when 
considering a sequence of initial conditions for which $|E_0|^2_{L^2}$ approaches $|R|^2_{L^2}
\approx 1.86$ from above. In this limit, $S$ 
is close to $R$, and $N$ close to $-R^2$.

\begin{figure}[t]
\centerline{
\includegraphics[height=7.5cm,width=8cm]{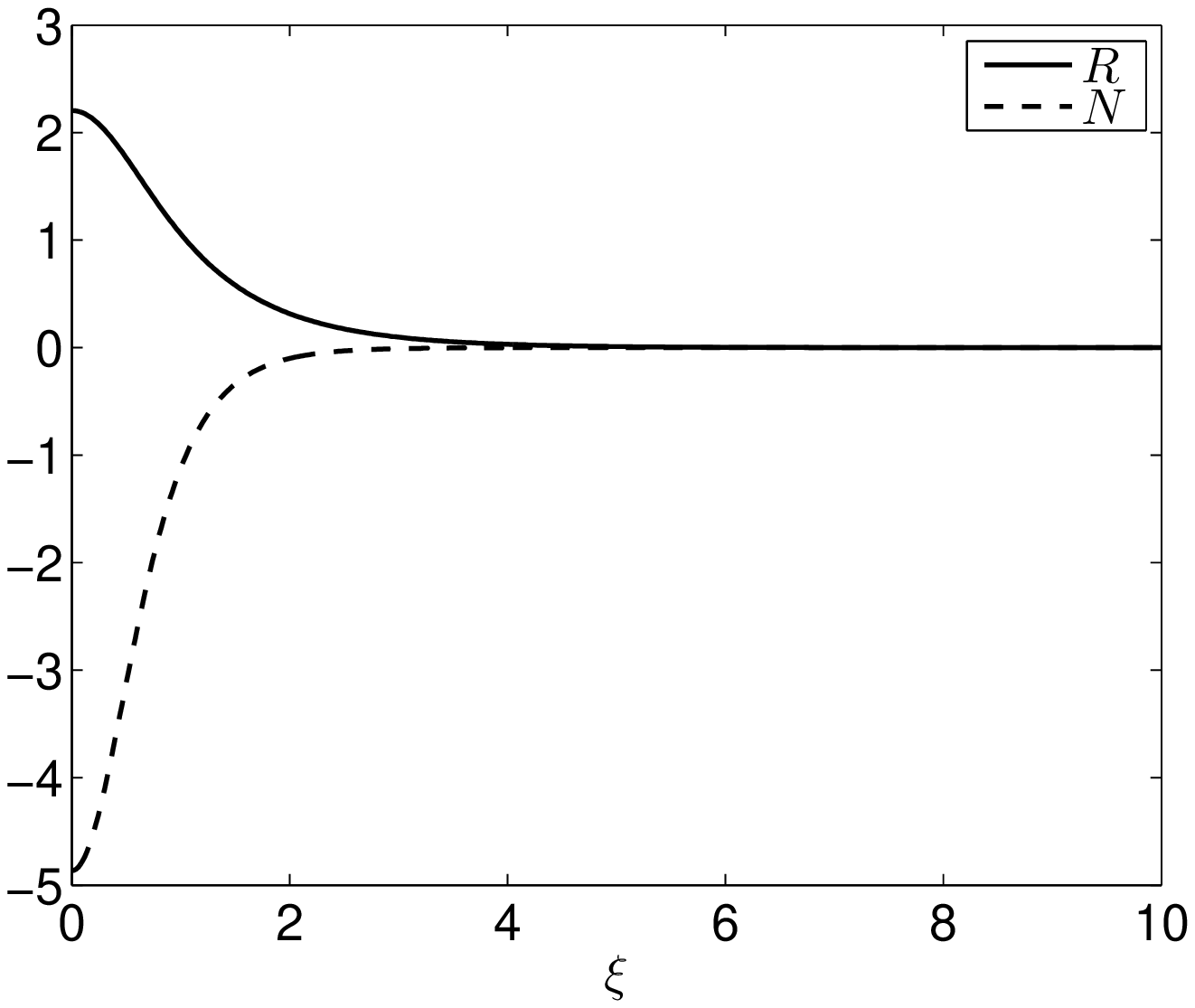}
}
\caption{Profiles of the positive solution (ground state) of the equation $\Delta R -R +R^3=0$
and of the density $N= -R^2$.}
\label{groundstate2}
\end{figure}

We now consider the effect of a small but non zero $\Gamma$, while keeping the assumption 
that the initial conditions are such that  $|E_0|^2_{L^2}$ is slightly in excess of  $|R|^2_{L^2}$. 
In this regime, the quantities  $a$ and $b$ are supposed to remain small but are allowed
to be time dependent. Furthermore, since the profiles of the various fields
remain
close
to those of the  blowing-up regime, we are led to expand them
in terms of the small parameters $a$, $b$ and $\Gamma/\lambda^2$ in the form
\begin{eqnarray}
&&S = R + a^2 \sigma_1 + b \sigma_2 + \frac{\Gamma}{\lambda^2} \sigma_3 +\cdots\\
&&N = -R^2 + a^2 \nu_1 + b \nu_2 + \frac{\Gamma}{\lambda^2} \nu_3 +\cdots\\
&&V = a \upsilon_1 + \cdots
\end{eqnarray}
For sake of completeness, the leading order profiles $R$ and $-R^2$
are  plotted in Fig \ref{groundstate2}.

At the center of symmetry, $\partial_\xi \sigma_i$,  $\partial_\xi \nu_i$ and $\upsilon_i$ vanish. All the 
functions also decay at infinity.
This leads to the sequence of equations
\begin{eqnarray}
&&\Delta \sigma_1 - \sigma_1 + R^2 \sigma_1 - R \nu_1 = 0\\
&& -\Delta \nu_1 - 2 \Delta (R\sigma_1)= {\mathcal L}(R^2)\\
&& \xi^{-1}\partial_\xi(\xi \upsilon_1)= 2R^2 + \xi (R^2)_\xi\\
&&\nonumber\\
&&\Delta \sigma_2 - \sigma_2 + R^2 \sigma_2 - R \nu_2 = -\frac{\xi^2}{4}R\\
&& -\Delta \nu_2 - 2 \Delta (R\sigma_2)= 0\\
&&\nonumber \\
&&\Delta \sigma_3 - \sigma_3 + R^2 \sigma_3 - R \nu_3 = \Delta^2R\\
&& -\Delta \nu_3 - 2 \Delta (R\sigma_3)= \Delta^2(R^2)
\end{eqnarray}
Since the kernel of the operator $\Delta - 1 + 3 R^2$ is reduced to the null function 
under a radial symmetry assumption,
all the above equations are solvable. We now write that  the 
time evolution of the fields through the variations of the functions  $a$, $b$ and $\lambda$
are constrained by the conservation of the plasmon number and of the Hamiltonian
that can both be estimated within  the above perturbative expansion \cite{Malkin93}. 

For the plasmon number, one has (up to a  $2\pi$ angular integration factor that we systematically omit)
\begin{eqnarray}
&&{\mathcal N} \equiv  \int |E|^2 r dr =\nonumber \\
&&\quad   \int \Big ( R^2  +2 a^2 R \sigma_1 + 2 b R\sigma_2 
+ 2 \frac{\Gamma}{\lambda^2} R \sigma_3  \Big) \xi  d\xi + \cdots.\nonumber \\
\label{plasmon}
&&
\end{eqnarray}
The Hamiltonian reads
\begin{eqnarray}
&&{\mathcal H} = \int \Big ( |\mbox{\boldmath $\nabla$}E|^2 + n |E|^2 + \frac{1}{2} n^2 + \frac{1}{2} V^2
\nonumber \\
&&\quad + \Gamma |\Delta E|^2 + \frac{\Gamma}{2} |\mbox{\boldmath $\nabla$} n|^2 \Big ) r dr.
\label{hamil}
\end{eqnarray}
It involves  
$\displaystyle|{\mbox{\boldmath $\nabla$} E|^2 =
\lambda^{-4} [|\mbox{\boldmath $\nabla$}S|^2 + (a^2\lambda^2/4)  \xi^2 S^2]}$.
Substituting the expansion of the various fields and using  that, as already mentioned,  the Hamiltonian associated 
with the ground state solution $(R, -R^2,0)$ vanishes, one gets
\begin{eqnarray}
&&\lambda^2 {\mathcal H} =  2 a^2 \int (\mbox{\boldmath $\nabla$}R \cdot \mbox{\boldmath $\nabla$} \sigma_1 
-R^3 \sigma_1) \xi d\xi \nonumber\\
&&\qquad + 2b \int (\mbox{\boldmath $\nabla$}R \cdot \mbox{\boldmath $\nabla$} \sigma_2 
-R^3 \sigma_2) \xi d\xi\nonumber\\
&&\qquad + \frac{2\Gamma^2}{\lambda^2} \int \mbox{\boldmath $\nabla$}R \cdot \mbox{\boldmath $\nabla$} \sigma_3 
-R^3 \sigma_3) \xi d\xi \nonumber\\
&&\qquad +\frac{a^2\lambda^2}{4} \int \xi^2R^2 \xi d\xi + \frac{a^2}{2} \int\upsilon_1^2  \xi d\xi\nonumber\\
&&\qquad +  \frac{\Gamma}{\lambda^2}\int|\Delta R|^2 \xi d\xi +\frac{\Gamma}{2\lambda^2}
\int \big|\mbox{\boldmath $\nabla$}(R)^2\big|^2 \xi d\xi + \cdots \nonumber\\
\end{eqnarray}
Performing integration by parts in the integrals involving the $\sigma_i$'s and using that 
$R$ obeys $\Delta R - R + R^3 =0$, one gets 
%(with ${\cal H} < 0$)
\begin{eqnarray}
&&\lambda^2  {\mathcal H} = -\int \Big (2a^2 R\sigma_1 + 2 b R \sigma_2 +  
\frac{2\Gamma}{\lambda^2}R \sigma_3 \Big) \xi d\xi \nonumber\\
&&+\frac{a^2\lambda^2}{4} \int \xi^2 R^2 \xi d\xi + \frac{a^2}{2} \int |\upsilon_1|^2 \xi d\xi \nonumber\\
&& +  \frac{\Gamma}{\lambda^2}\int\Big (|\Delta R|^2  
+\frac{1}{2} \big|\mbox{\boldmath $\nabla$}(R^2)\big|^2  \Big)
\xi d\xi + \cdots \label{Hexp}
\end{eqnarray}
Using (\ref{plasmon}), 
the first integral  in the right-hand side of  eq. (\ref{Hexp}) is the difference ${\widetilde {\mathcal N}}=
{\mathcal N} -\int R^2 \xi d\xi$ between the actual plasmon number and its critical value for collapse when $\Gamma = 0$. 
Furthermore,   $a^2 = \lambda_t^2$. One thus gets the following dynamical equation for the scaling factor $\lambda$
\begin{equation}
\lambda_t^2 = \frac{1}{m_1 + m_2 \lambda^2} \Big ( {\mathcal H} \lambda^2+ {\widetilde {\mathcal N}}- 
\frac{\Gamma m_3}{\lambda^2}  \Big),
\end{equation}
with $m_1 = (1/2) \int \upsilon^2 \xi d\xi\approx 0.727$, $m_2 =(1/4) \int \xi^2 R^2\xi d\xi \approx 0.553$ 
and $m_3 = \int \Big (|\Delta R|^2+ (1/2) \big|\mbox{\boldmath $\nabla$}(R^2)\big|^2 \Big ) \xi d\xi \approx 10.785$ and  ${\cal H} < 0$.
As in \cite{HS09}, the problem can be viewed as
the motion of a particle in a potential that diverges positively as $\lambda \to 0$
and has a finite positive limit as $\lambda \to +\infty$. If $\Gamma < {\widetilde {\cal N}}^2/4|{\cal H}|m_3$,
this potential has a negative minimum and 
$\lambda(t)$ oscillates between  strictly positive minimal and maximal  values.

\begin{figure}[t]
%\centerline{
%\includegraphics[height=7.5cm,width=8cm]{figure2a.ps}
%}
%\centerline{
%\includegraphics[height=7.5cm,width=8cm]{figure2b.ps}
%}
\centerline{
\includegraphics[height=16cm,width=8cm]{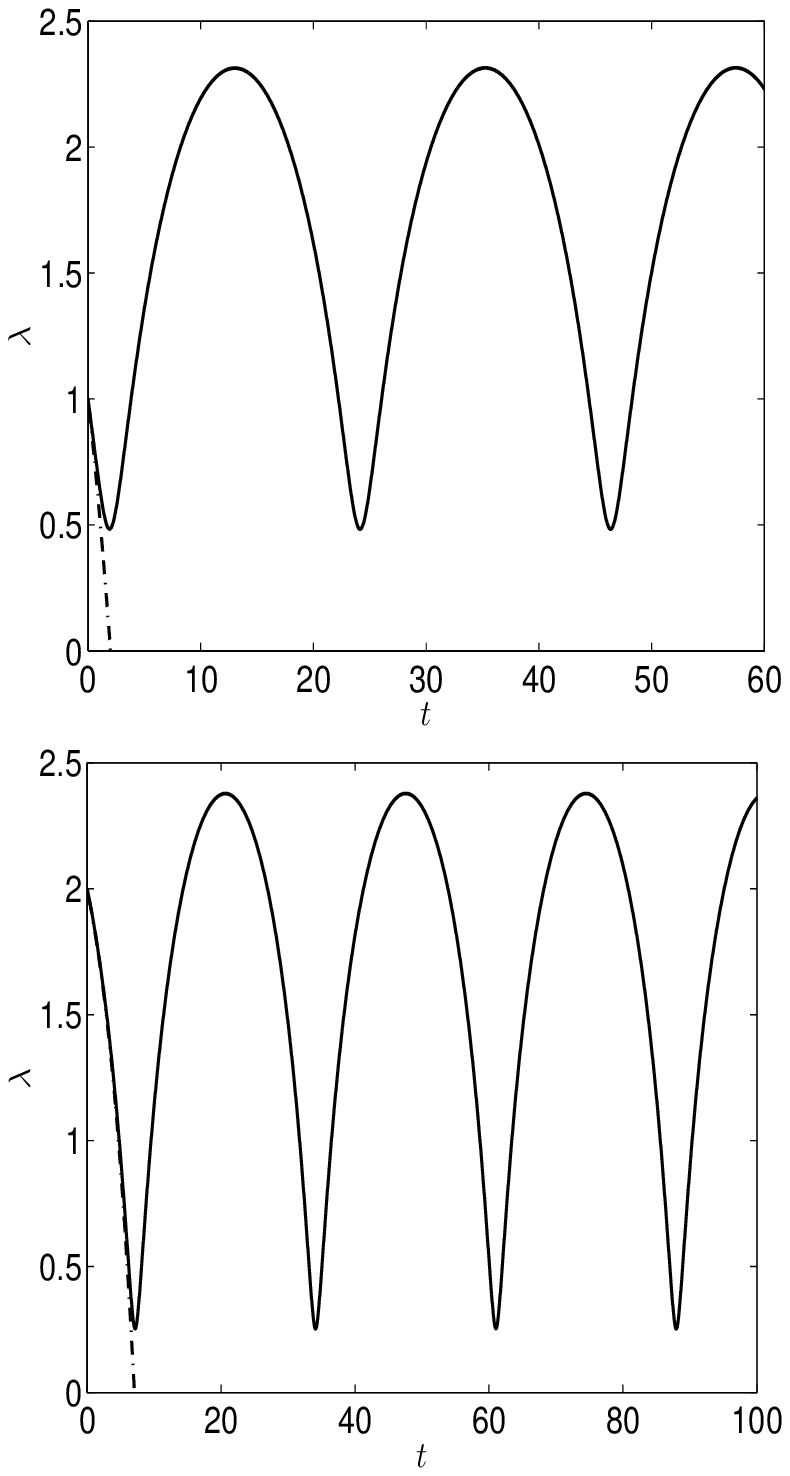}
}
\caption{Time variation of the scaling factor $\lambda$ (solid lines) for 
$\Gamma = 5. \ 10^{-3}$, ${\mathcal H}= -0.0430  $ and ${\widetilde {\mathcal N}}= 0.240$ (top) and  
$\Gamma = 10^{-3}$, ${\mathcal H}= - 0.0295 $ and ${\widetilde {\mathcal N}}= 0.168$  (bottom).
For comparison, the corresponding evolutions in the case $\Gamma = 0$ (dashed lines) are superimposed.}
\label{osc}
\end{figure}

For an explicit numerical resolution, it is more convenient to consider $y= \lambda^2$
that obeys 
\begin{equation}
y_t^2 =\frac{4}{m_1 + m_2 y}\Big ( {\mathcal H}y^2 +{\widetilde {\mathcal N}}y -\Gamma m_3 \Big ) \label{eqy2}
\end{equation}
Under the condition $\displaystyle{\Gamma < 
{\widetilde {\mathcal N}^2}/\big (4\big |{\mathcal H}\big |m_3\big )}$, there exists positive $y_m$ and $y_M$  solutions of 
${\mathcal H} y^2+ {\widetilde {\mathcal N}}y - {\Gamma m_3} $ such that $y$ oscillates 
between $y_m= \displaystyle{[{\widetilde {\cal N}}-({\widetilde{\cal N}}^2
-4{|\cal H|}\Gamma m_3)^{1/2}]/2{|\cal H|}}$ 
and $y_M= \displaystyle{[{\widetilde {\cal N}}+({\widetilde{\cal N}}^2
-4{|\cal H|}\Gamma m_3)^{1/2}]/2{|\cal H|}} $, 
consistent with the global existence demonstrated in Section III. Equation (\ref{eqy2})
is equivalently rewritten
\begin{equation}
y_{tt}= \partial_y \Big [\frac{2}{m_1 + m_2 y}\Big ( {\mathcal H}y^2 +{\widetilde {\mathcal N}}y -\Gamma m_3 \Big )\Big ].
\end{equation}
For a given $y(0)$, one should prescribe 
\begin{equation}
y_t(0)= - 2\left [\frac{ {\mathcal H}y^2(0) +{\widetilde {\mathcal N}}y(0) -\Gamma m_3 }{m_1 + m_2 y(0)}\right]^{1/2}.
\end{equation}

Figure \ref{osc}  shows the time evolution of the scaling factor $\lambda(t)$ for  $\Gamma = 5. \ 10^{-3}$, ${\mathcal H}=  -0.0430$ and ${\widetilde {\mathcal N}}= 0.240$ (top) and  
$\Gamma = 10^{-3}$, ${\mathcal H}= - 0.0295$ and ${\widetilde {\mathcal N}}= 0.168$  (bottom), 
corresponding to the initial  
conditions  $E_0 = 2.85 \,e^{-r^2}$ and 
$E_0 = 2.90 \,e^{-r^2}$  respectively, with $n_0= -|E_0|^2$, $V_0=0$. 
For comparison, we superimposed the corresponding evolution when $\Gamma=0$. In this case, 
$\lambda$ reaches zero in a finite 
time, with  $\lambda_t^2 ={\widetilde N}/m_1=a_0^2$, consistent with the  scaling   $\lambda \propto (t_*-t)$
of the self-similar blowing-up solutions of the classical Zakharov equations. 
Note that, while near threshold  the
two-dimensional scalar model predicts the  same leading- order profile for the pump wave  as the  nonlinear Schr\"odinger equation 
resulting from the subsonic approximation (slaved acoustic waves), it does not lead to the same scaling law.

\section{The two-dimensional electrostatic model}
Let us now return to the electrostatic model, first in the
case $\Gamma=0$. From eqs. (\ref{zakhpot1})-(\ref{zakhpot2}), we perform the rescaling $\displaystyle{
\nabla\psi = (1/\lambda) {\mathbf S}({\mbox{\boldmath $\xi$}, \tau) e^{i\tau -ia \lambda {|\mbox{\boldmath $\xi$}|^2}/4}}}$,
$\displaystyle{n = (1/\lambda^2) N(\mbox{\boldmath $\xi$}, \tau)}$ and  
$\displaystyle{{\mathbf v} = (1/\lambda^2) {\mathbf V}(\mbox{\boldmath $\xi$}, \tau)}$.
After neglecting as previously the $\tau$-derivatives of the rescaled functions, we get
\begin{eqnarray}
&&\Delta {\mathbf S} -  {\mathbf S} + b \frac{|\mbox{\boldmath $\xi$}|^2}{4}
{\mathbf S} =\nonumber \\
&&\qquad e^{ia \lambda {|\mbox{\boldmath $\xi$}|^2}/4} \Delta^{-1}  \mbox{\boldmath $\nabla$} 
\Big (\mbox{\boldmath $\nabla$}\cdot(N {\mathbf S} \,e^{-ia \lambda {|\mbox{\boldmath $\xi$}|^2}/4} )\Big) \label {bfS}\\
&& a^2 {\cal L} (N) - \Delta N = \Delta |{\mathbf S}|^2. 
\end{eqnarray}
%To leading order, we have to solve 
%\begin{equation}
%\Delta {\mathbf R} -  {\mathbf R} +  \Delta^{-1}  \mbox{\boldmath $\nabla$} 
% (\mbox{\boldmath $\nabla$}\cdot \Big (|{\mathbf R}|^2 {\mathbf R})\Big )=0 \label {bfS}\\
%\end{equation}
%but the presence of the exponential introduces an additional equation compared with the
%scalar case, of the form
%\begin{equation}
%\Delta {\mathbf s} -  {\mathbf s} = -i a \lambda}{\mathbf F}
%\end{equation}
%with 
%\begin{equation}
%\frac{|\mbox{\boldmath $\xi$}|^2}{4} \Delta^{-1}  \mbox{\boldmath $\nabla$} 
% (\mbox{\boldmath $\nabla$}\cdot (|{\mathbf R}|^2{\mathbf R})
%+\Delta^{-1}  \mbox{\boldmath $\nabla$} 
% \Big [\mbox{\boldmath $\nabla$}\cdot \Big (\frac{|\mbox{\boldmath $\xi$}|^2}{4}
%|{\mathbf R}|^2 {\mathbf R}\Big ) \Big]\\
%\end{equation}

The phase factors in eq. (\ref{bfS}) introduce a serious difficulty in the sense that 
their expansions  lead to an additional contribution in  the perturbative
calculation that is  not necessarily associated with a well-posed problem, 
restricting de facto the present analysis to isotropic
solutions for which the operator $\Delta^{-1} \, {\rm grad} \, {\rm div}$ 
reduces to the identity and  the phase factors 
cancel out. Denoting by $S$ the radial component of ${\mathbf S}$, one can then perform an 
analysis  similar to that of the scalar model.

\begin{figure}[t]
\centerline{
\includegraphics[height=7.5cm,width=8cm]{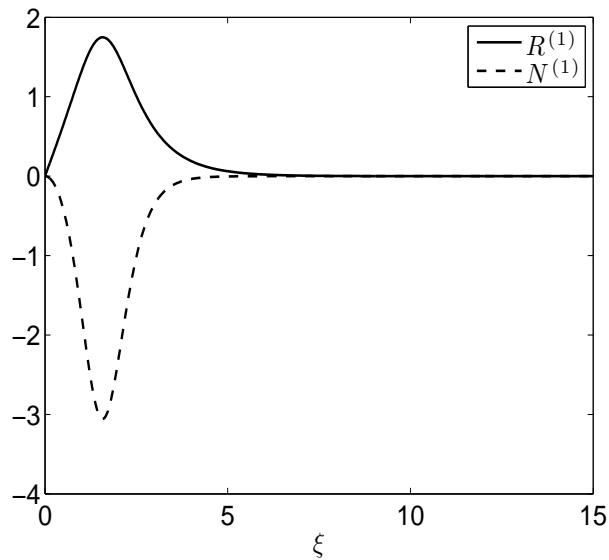}
}
\caption{Profiles of the positive solution (ground state) of the equation 
$\Delta^{(1)} R^{(1)} -R^{(1)} +{R^{(1)}}^3=0$ and of the density $N^{(1)}= -{R^{(1)}}^2$.}
\label{vortex1}
\end{figure}

The only difference between the isotropic electrostatic equations and
the scalar  model is the replacement of the scalar Laplacian in the equation for the
electric field  by the radial component ${ \Delta^{(1)}}$ of the vectorial Laplacian 
(which implies in particular that the electric field now vanishes at the center of symmetry). Denoting
by $R^{(1)}$ the positive solution (fig. \ref{vortex1}) of 
\begin{equation}
\Delta^{(1)} R^{(1)} - R^{(1)} + {R^{(1)}}^3 =0
\end{equation}
with $R^{(1)}(0)= R^{(1)}(\infty) = 0$ (see \cite{FG08} for a review on this equation and its extensions
that arise in the context of vortex solutions of the nonlinear Schr\"odinger equation),
the rescaled density and velocity satisfy the same equation as in the scalar model,
except that $R$ is replaced by $R^{(1)}$. We thus recover eq. (\ref{eqy2}) with coefficients now given by 
$m_1 \approx 24.42$, $m_2  \approx 8.14$ and $m_3  \approx 24.94$.
Furthermore $\int {R^{(1)}}^2\xi d\xi \approx 7.69$.
As  consequence, the scaling coefficient displays the same oscillatory behavior as in the scalar model.

\section{Difficulties in three dimensions}
The three-dimensional problem displays specific difficulties even in the context of the scalar model. 
Indeed, although solutions of the classical scalar model ($\Gamma =0$)
with negative Hamiltonian appear to blow up in a self-similar way  as in two dimensions,
the perturbation analysis for $\Gamma >0$ is not  straightforward
in three dimensions. 

Proceeding as in two dimensions (but with different rescalings), we  are looking for solutions of the form
$E(r,t) = \lambda^{-3/2} S(\xi, \tau)e^{i(\tau - a \lambda^{1/2} \xi^2/4)}$, 
$n(r,t) = \lambda^{-2} N(\xi, \tau)$,
$v(r,t) = \lambda^{-5/2} V(\xi, \tau)$, 
$ n_t(r,t) = \lambda^{-7/2} W(\xi, \tau)$, 
with  $\xi = r/\lambda$ with a profile depending only weakly on the rescaled time $\tau$. This leads to 
\begin{eqnarray}
&&\Delta S - S - NS + b\frac{\xi^2}{4} S \nonumber \\
&& \qquad\qquad - \frac{\Gamma}{\lambda^2} e^{ia\lambda\xi^2/4}\Delta (Se^{-ia\lambda\xi^2/4}) = 0 \label{S3d}\\
&&a( \frac{7}{2}W + \xi \partial_\xi W ) - \lambda \Delta N = \Delta |S|^2 - \frac{\Gamma}{\lambda} \Delta^2 N\\
&& \frac{1}{\xi^2}\partial_\xi (\xi^2 \partial_\xi V) = -a (2 N + \xi \partial_\xi N) = -W\label{V3d}\\
&& a ( \frac{5}{2}V + \xi \partial_\xi V )- \lambda \partial_\xi N =
- \partial_\xi| S|^2 
+ \frac{\Gamma}{\lambda} \partial_\xi \Delta N \label{Vt3d}
\end{eqnarray}
where $a=-\lambda_t\lambda^{1/2} = - \lambda_\tau \lambda^{-3/2}$ and 
$b=(a\lambda^{1/2})_\tau + (a\lambda^{1/2})^2$. From eqs. (\ref{V3d}) and  (\ref{Vt3d}), we have 
\begin{equation}
a^2 {\mathcal L}(N) - \lambda\Delta N = \Delta |S|^2 - \frac{\Gamma}{\lambda}\Delta^2 N\label{N3d} \\
\end{equation}
with now  ${\cal L} = \xi^2 N_{\xi\xi} + (13/2) \xi N_\xi + 7 N$.

In the classical regime ($\Gamma = 0$), $a$ has a finite (positive) limit $a_0$ as the collapse time is approached,
corresponding to a scaling factor varying like $(t_\star-t)^{2/3}$. In this regime, the term $\lambda \Delta N$
is negligible in eq. (\ref{N3d}) (supersonic regime) 
and $S$ scales like $a$. However, as shown below, this ordering 
eventually breaks down  in the presence of quantum effects. Indeed, as the collapse is arrested, $\lambda_t$
vanishes and so does $a$. This indicates that a systematic modulational theory analogous to what we developed 
in two dimensions, is not possible in three dimensions. In this context, we resort to limit our address of the
problem to a semi-phenomenological approach based on the extension of asymptotic expansions outside their 
range of strict validity, with the hope to capture qualitative properties of the global dynamics.

For this purpose, we assume that $\Gamma$ is small enough for the quantum effects to start acting 
only after  the 
system has reached the classical blowing up regime where $a(t)$ is closed to its limit $a_0$. In this asymptotic regime, one  expects that there exists a period of time during which 
$a$ remains sufficiently close to $a_0$ to allow a perturbative calculation. We are thus led to expand
\begin{eqnarray}
&& S = a_0 (S_0 + b S_1 + \frac{\lambda}{a_0^2} S_2 + \frac{\Gamma}{\lambda^2} S_3 +
\frac{\Gamma}{\lambda a_0^2} S_4 + \frac{a^2-a_0^2}{a_0^2} S_5 )\nonumber \\
&& N =  N_0 + b N_1 + \frac{\lambda}{a_0^2} N_2 + \frac{\Gamma}{\lambda^2} N_3 +\frac{\Gamma}{\lambda a_0^2} N_4+
\frac{a^2-a_0^2}{a_0^2}N_5 \nonumber\\
&&V= a_0(V_0 + b V_1 + \frac{\lambda}{a_0^2} V_2 + \frac{\Gamma}{\lambda^2} V_3 +  \frac{\Gamma}{\lambda a_0^2} V_4+
\frac{a^2-a_0^2}{a_0^2}V_5 \nonumber \\
&&\qquad + (a-a_0) V_0\nonumber).
\end{eqnarray}

The lowest order terms $(S_0,N_0, V_0)$ are solutions of 
\begin{eqnarray}
&&\Delta S_0 - S_0 - N_0 S_0 =  0 \label{s0}\\
&& {\mathcal L}(N_0) =  \Delta (S_0^2)  \label{n0}\\
&& -\xi^{-2} \partial_{\xi}  (\xi^2 V_0)= 2N_0 + \xi \partial_\xi N_0.\label{v0}\\
&& \frac{5}{2} V_0 + \xi {V_0}_\xi = - \partial_\xi (|S_0|^2) \label{vt0}
\end{eqnarray}
Their profiles (not shown) are  qualitatively very similar to those in two dimensions \cite{LPS92}.
The corrections terms $(S_i,N_i, V_i)$ $(i=1,5)$ satisfy the systems
\begin{eqnarray}
&&\Delta S_i - S_i - N_0 S_i - S_0N_i = F_i \label{si}\\
&& {\mathcal L}(N_i) -  2 \Delta (S_0S_i) =G_i \label{n1} \\
&& \xi^{-2} \partial_{\xi}  (\xi^2 V_i)= H_i,\label{vi}
\end{eqnarray}
with 
$F_1 =-\frac{\xi^2}{4} S_0$, $G_1= 0$, 
$F_2=0$, $G_2= \Delta N_0$, 
$F_3= \Delta^2 S_0$, $G_3= 0$, 
$F_4= 0$, $G_4= -\Delta^2N_0$,
$F_5=0$, $G_5 = - {\mathcal L}(N_0)$,
and  $H_i= -a_0(2 N_i + \xi \partial_\xi N_i)$.
We find that $S_5=\frac{1}{2} S_0$ and $N_5=V_5=0$. The first and last term in the expansion
of $S$ thus combine giving $\frac{a^2+a_0^2}{2a_0} S_0$.
Proceeding as in two dimensions, we
substitute the expansions of $S,N,V$ into the plasmon number and the
Hamiltonian.  We get 
\begin{equation}
{\mathcal N} = \alpha _0 a^2  +2\alpha _1 a_0^2 b  + 2 \alpha _2\lambda 
+ 2 \alpha_3 \Gamma \frac{a_0^2}{\lambda^2}  +2 \alpha_4  \frac{\Gamma}{\lambda}   + \cdots.
\label{plasmon2}
\end{equation}
where the coefficients $\alpha_i \equiv \int S_0 S_i \xi^2 d\xi$  are numerically estimated as
$\alpha_0 \approx 1.99$, $\alpha_1 \approx 1.94$,
$\alpha_2 \approx 0.35$,  $\alpha_3 \approx 51.36$,  $\alpha_4 \approx 34.44$.
On the other hand,
\begin{equation}
\lambda^2 {\mathcal H} =  \beta_0 a_0^2 +  \beta_1  a_0^2a^2\lambda +  \beta_2 a_0^2b+ \beta_3 \lambda 
+ \beta_4 \Gamma  \frac{a_0^2}{\lambda^2} +  \beta_5 \frac{\Gamma}{\lambda} + \cdots
\label{hamil3}
\end{equation}
One first proves that 
\begin{equation}
\beta_0 \equiv \int\Big(|\mbox{\boldmath $\nabla$}S_0|^2 + n |S_0|^2 + \frac{1}{2} V_0^2 \Big)\xi^2  d\xi =0. \label{inv3d}
\end{equation}
For this purpose, we multiply eq. (\ref{s0}) by  $S_0$ and by $ \xi (S_0)_\xi$ respectively,
and integrate in space
the resulting equations to get
\begin{eqnarray}
\int ( |\mbox{\boldmath $\nabla$}S_0|^2 + N_0 S_0^2 +  S_0^2 )\xi^2  d\xi = &0&,
\label{A1} \\
\int ( |\mbox{\boldmath $\nabla$}S_0|^2 + 3 N_0 S_0^2 + 3 S_0^2
+ \xi  {N_0}_\xi  S_0^2  )\xi^2  d\xi =&0&.
\label{A2}
\end{eqnarray}
Subtracting these two equalities and  using (\ref{v0}),
we have 
\begin{equation}
\int [  2 S_0^2 + V_0  \partial_\xi(S_0^2)]  \xi^2  d\xi =0.
\label{A4}
\end{equation}
Multiplying (\ref{vt0}) by $V_0$ and integrating  in space, we get
\begin{equation}
\int V_0^2 \xi^2  d\xi =- \int (S_0^2)_\xi V_0\xi^2  d\xi .
\label{A5}
\end{equation}
Combining (\ref{A4}) and (\ref{A5}) gives $\int S_0^2 \xi^2d\xi = \frac{1}{2}
\int V_0^2 \xi^2  d\xi$, which after substitution in 
(\ref{A1}) leads to  (\ref{inv3d}).

Using eqs.  (\ref{s0})-(\ref{vi}),
the constants $\beta_i$ for $i=1,\cdots 5$  are given by
% (omitting  the volume element $\xi^2d\xi$)
%\begin{eqnarray}
%&& \beta_1 = \int \frac{\xi^2}{4} S_0^2 \approx ...  \\
%&& \beta_2 = \int\Big(2 \nabla S_0\cdot\nabla S_1 + 2N_0S_0S_1 + N_1 S_0^2 + V_0V_1\Big) \nonumber\\
%&&\quad =  \int\Big(\frac{\xi^2}{4} S_0^2 -2 S_0S_1 +V_0V_1 \Big)\approx ...  \\
%&&\beta _3 = \int\Big(2 \nabla S_0\cdot\nabla S_2 + 2N_0S_0S_2 + N_2 S_0^2 + V_0V_2\nonumber\\
%&& \quad +\frac{1}{2} N_0^2 \Big)
% =  \int\Big( -2 S_0S_2 +V_0V_2 +\frac{1}{2} N_0^2  \Big) \approx ... \\
%&& \beta_4 =\int\Big(2 \nabla S_0\cdot\nabla S_3+ 2N_0S_0S_2 + N_3 S_0^2 + V_0V_3
%\nonumber\\
%&&\quad  + (\Delta S_0)^2\Big)= \int\Big( -2 S_0S_3 +V_0V_3 \Big)\approx ...  \\
%&&\beta_5 =  \int \Big( 2 \nabla S_0\cdot \nabla S_4 + 2N_0S_0S_4 + N_4 S_0^2+ V_0V_4\nonumber\\
%&& \quad + \frac{1}{2}|\nabla N_0|^2 \Big)
% = \int \Big(-2 S_0S_4 +V_0V_4 +\frac{1}{2}|\nabla N_0|^2 \Big) \nonumber\\
%&&\qquad \approx ... 
%\end{eqnarray}
%\begin{eqnarray}
%&& \beta_1 = \int \frac{\xi^2}{4} S_0^2 \approx ...  \\
%&& \beta_2 = \int(\frac{\xi^2}{4} S_0^2 -2 S_0S_1 +V_0V_1 )\xi^2d\xi \approx ...  \\
%&&\beta _3 =  \int ( -2 S_0S_2 +V_0V_2 +\frac{1}{2} N_0^2 )\xi^2d\xi \approx ... \\
%&& \beta_4 =  \int( -2 S_0S_3 +V_0V_3 )\xi^2d\xi \approx ...  \\
%&& \beta_5  = \int (-2 S_0S_4 +V_0V_4 +\frac{1}{2}|\mbox{\boldmath $\nabla$} N_0|^2 )\xi^2d\xi
%\end{eqnarray}
$ \beta_1 = \int \frac{\xi^2}{4} S_0^2 \approx 1.17  $, 
$\beta_2 = \int(\frac{\xi^2}{4} S_0^2 -2 S_0S_1 +V_0V_1 )\xi^2d\xi \approx -1.17 $,
$\beta _3 =  \int ( -2 S_0S_2 +V_0V_2 +\frac{1}{2} N_0^2 )\xi^2d\xi \approx 3.10$,
$ \beta_4 =  \int( -2 S_0S_3 +V_0V_3 )\xi^2d\xi \approx -5.56 $,
$\beta_5  = \int (-2 S_0S_4 +V_0V_4 +\frac{1}{2}|\mbox{\boldmath $\nabla$} N_0|^2 )\xi^2d\xi
\approx -2.81$.

We now eliminate the term proportional to $a^2b$ in ${\mathcal H} $ using ${\mathcal N}$ and find
\begin{eqnarray}
&&\lambda^2 {\mathcal H} =   m_1  a_0^2 a^2\lambda + m_2 \lambda 
+m_3\Gamma   \frac{a_0^2}{\lambda^2} +m_4 a^2 \nonumber\\
&&\quad +m_5\frac{\Gamma}{\lambda} + m_6{\cal N}.
\label{hamil4}
\end{eqnarray}
The new constants $m_i$ are defined as 
$m_1 = \beta_1\approx 1.17$, $m_2= \beta_3 - \alpha_2 \beta_2 /\alpha_1\approx 3.31$, 
$m_3= \beta_4 - \alpha_3\beta_2 /\alpha_1\approx 25.41$,
$m_4= -\alpha_0\beta_2 /(2\alpha_1)\approx 0.60$, $m_5= \beta_5 - \alpha_4\beta_2 /\alpha_1\approx 
17.95$, 
$m_6 = \beta_2 /(2\alpha_1)\approx -0.30$. 

Inserting that $a=-\lambda_t\lambda^{1/2}$,
the effective ODE satisfied by $\lambda$ is
\begin{equation}
\lambda_t^2 = -\frac{ m_2\lambda^3+ m_3\Gamma a_0^2+ m_5\Gamma\lambda+ m_6{\mathcal N}\lambda^2- 
{\mathcal H}\lambda^4}{\lambda^2( m_1 a_0^2 \lambda^2+ m_4\lambda  )}
\label{ode2}
\end{equation}
In order to fix the parameter $a_0$, we note that for $\Gamma = 0$, in the blowing up regime where $\lambda \to 0$,
one has 
\begin{equation}
a_0^2 = \lim_{\lambda\to 0} a^2 =  \lim_{\lambda\to 0} \lambda_t^2 \lambda = -\frac{ m_6{\mathcal N}}{ m_4},
\end{equation}
or using the definitions of $m_4$ and $m_6$,
\begin{equation}
a_0^2 = \frac{\mathcal N}{\int S_0^2 \xi^2 d\xi}.
\end{equation}

\begin{figure}[t]
%\centerline{
%\includegraphics[height=7.5cm,width=8cm]{figure4a.ps}
%}
%\centerline{
%\includegraphics[height=7.5cm,width=8cm]{figure4b.ps}
%}
\centerline{
\includegraphics[height=16cm,width=8cm]{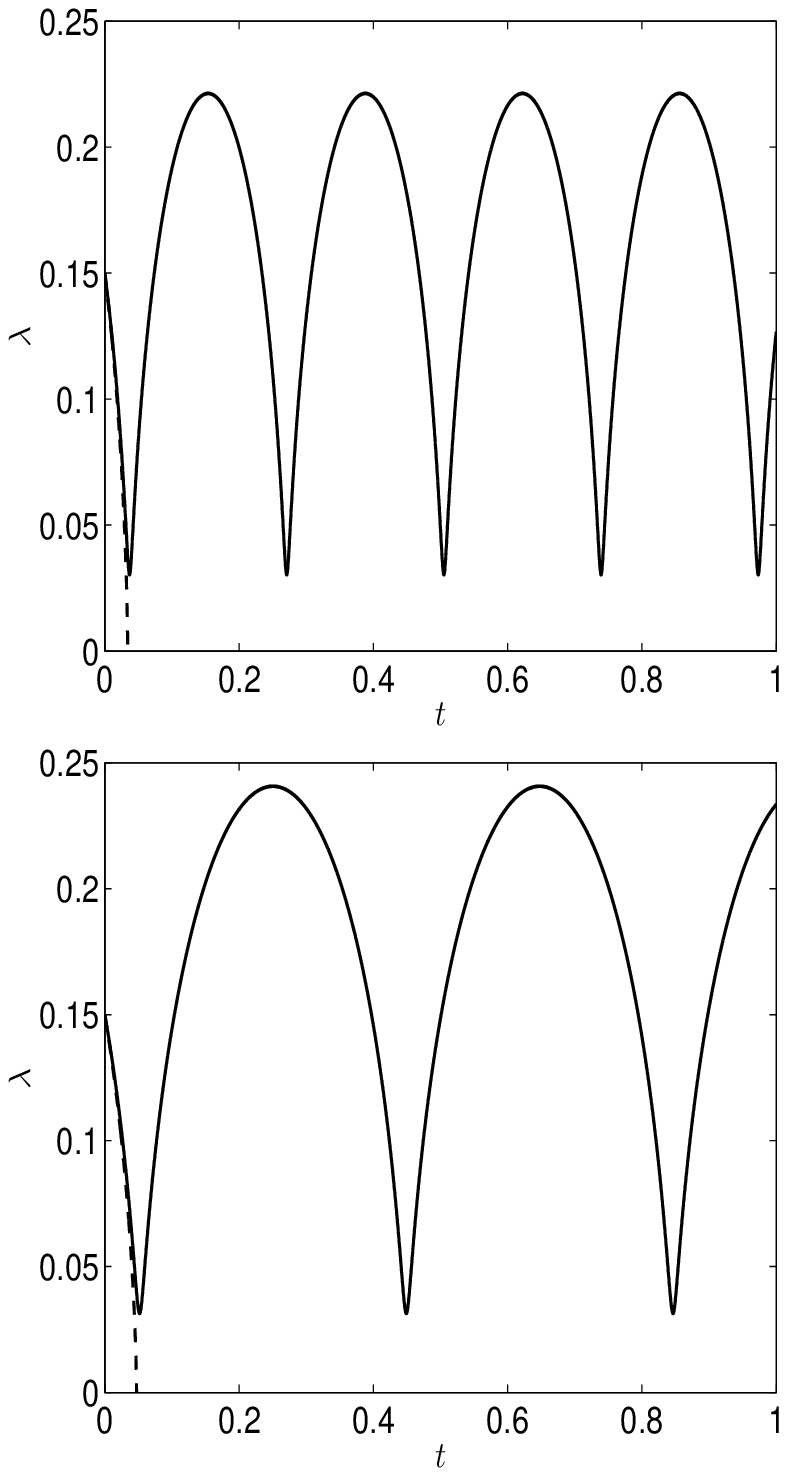}
}
\caption{Time variation of the scaling factor $\lambda$ 
for the three-dimensional phenomenological model, when  $ {\mathcal N} =5.64 $ and ${\mathcal H} = -18.97 $,
with $\Gamma = 2. \ 10^{-5}$ (top) and
$ {\mathcal N} =2.76 $ and ${\mathcal H} = -0.33$,
with the sale value of $\Gamma$ as above. For comparison, the corresponding evolutions in  the case $\Gamma =0$
(dashed lines) are  superimposed.}
\label{3d}
\end{figure}

Figure  4 (top)  shows the evolution of the scaling factor $\lambda(t)$
for  $ {\mathcal N} =5.64 $ and ${\mathcal H} = -18.97 $ (corresponding to an initial conditions 
$E_0 = c e^{-r^2} $, $n_0 =-|E_0|^2$, and $v_0=0$ with $c= 6$,
for the original Zakharov equations) and   $\Gamma = 2. \ 10^{-5}$.
For comparison, the evolution in the absence of quantum effects
is also displayed. As in two dimensions, we observe that quantum effects induce a periodic behavior.
In simulations with smaller $\Gamma$, the minimum of $\lambda$ is, as expected, getting smaller. 
When keeping the same value of $\Gamma$, one uses   ${\mathcal H} = -0.32 $ 
with ${\mathcal N} = 2.76$ by  
taking $c=4.2$, the maximum of $\lambda$ is slightly increased,  while the period
of the oscillation gets significantly longer (Fig. 4 bottom), but  the global behavior remains  very similar. 
Note however that the dynamics is much faster than in two dimensions. This effect is already visible on the singularity time
at $\Gamma=0$, a regime for which eq. (\ref{ode2}) is asymptotically exact. This suggests that the parameters we used in
three dimensions (even for the bottom panel of Fig. 4) correspond to a regime significantly distant from the
threshold conditions that in this case are not as easily characterized  as in two dimensions.
As stressed at the beginning of this section, the present description of the three-dimensional 
problem is however to be viewed as  heuristic, 
as the asymptotics clearly breaks down before the scaling factor $\lambda(t)$ reaches its minimum.
It nevertheless predicts a behavior consistent with the arrest of collapse (the proof given in Section I
is easily transposed to the scalar model). It also shows a  periodic dynamics whose origin is 
expected to be generic. It indeed  results from a competition between 
wave focusing that occurs  when $\lambda$ is not yet small enough for the quantum effects to act efficiently, and
the subsequent evolution that  takes place when the influence of the latter 
perturbations dominates the dynamics and leads to defocusing until the moment where, $\lambda$ becoming large enough,
their influence  becomes again subdominant, thus permitting an efficient self-focusing. 
Validation of the model would  require comparisons with direct simulations of the scalar model, 
an issue that is outside the scope of the present paper. Compared with the Rayleigh-Ritz method, the present approach
should provide a better description of the solution profile. It also more clearly points out the conditions of
applicability of modulation methods.

\section{Conclusion}
The influence of quantum effects on the Langmuir wave dynamics provides an interesting example of the action of an
additional dispersive effect on the phenomenon of wave collapse, although in realistic situations 
the Zakharov description is supposed to break down before quantum effects become relevant. 
Arrest of collapse was predicted \cite{HS09} in the adiabatic regime where the density is slaved to the wave amplitude, using 
a Rayleigh Ritz method.  Here, the result 
is rigorously established for the full quantum Zakharov equations, by combining the 
conservation of the plasmon number and of the Hamiltonian with estimates based on a Gagliardo-Nirenberg inequality.  
These invariances are also used to develop a systematic perturbative expansion in order to capture the influence 
of weak quantum effects for initial conditions slightly above the singularity threshold for the classical problem. 
Restricted for technical reasons to isotropic solutions (at least in the focusing region), 
the analysis is carried out in two space dimensions corresponding  to the critical dimension for collapse 
of the cubic nonlinear Schr\"odinger equation. 
The difficulty of extending the analysis to three dimensions points out the importance of the proximity of
criticality to satisfy the delicate balances involved in a systematic asymptotic theory. 
We thus resorted in this case to develop a semi-phenomenological approach.

\begin{acknowledgments}
We thank J. Colliander and G. Fibich for useful comments. 
This work was  partially supported by NSERC through grant number 46179-05.
\end{acknowledgments}

\bibliographystyle{apsrev}

\bibliography{quantum_zak}

\end{document}